\documentclass[12pt]{amsart}

\usepackage{amsmath,graphics,comment,eufrak}

  \includecomment{inlong} \excludecomment{inshort}   
 \includecomment{inlow} \excludecomment{inhigh}

\theoremstyle{plain}
\newtheorem*{theorem*}{Theorem}
\newtheorem*{lemma*} {Lemma}
\newtheorem*{corollary*} {Corollary}
\newtheorem*{proposition*} {Proposition}
\newtheorem{theorem}{Theorem}[section]

\newtheorem{corollary}[theorem]{Corollary}
\newtheorem{proposition}[theorem]{Proposition}

\theoremstyle{remark}

\newtheorem*{definition}{Definition}

\theoremstyle{definition}

\def \Z {\mathbf{Z}}
\def \C {\mathbf{C}}

\def\eps{\epsilon}

\def\s{\sigma}

\def\Q{\Bbb{Q}}
\def\id{\mbox{id}}

\def\sign{\mbox{sign}}
\def\Z{\Bbb{Z}}
\def\C{\Bbb{C}}

\def\N{\Bbb{N}}

\def\part{\partial}

\def\a{\alpha}

\def\mapstoleft{\leftarrow \hspace{0.2cm} \line(0,1){5.3} \hspace{-0.16cm}\times}
\def\bp{\begin{pmatrix}}

\def\sm{\setminus}
\def\ep{\end{pmatrix}}
\def\bn{\begin{enumerate}}

\def\en{\end{enumerate}}
\def\ba{\begin{array}}
\def\ea{\end{array}}

\def\s{\sigma}
\def\a{\alpha}

\def\ti{\tilde}

\def\fr12{\frac{1}{2}}

\def\Oplus{\bigoplus}
\def\im{\mbox{Im}}

\def\ker{\mbox{Ker}}
\def\tr{\mbox{tr}}
\def\t{\theta}

\def\hom{\mbox{Hom}}

\def\mapstoleft{\leftarrow \hspace{-0.14cm}\line(0,1){5.3}}

\def\cmtbf#1{} \def\cmt#1{}

\begin{document}

\title{Full signature invariants for $L_0(F(t))$}
\author{Stefan Friedl}
\date{\today}
\begin{abstract}
Let $F/\Q$ be a number field closed under complex conjugation.
Denote by ${L}_0(F(t))$ the Witt group of hermitian
forms over $F(t)$. 
  We find full invariants for detecting non--zero elements in ${L}_0(F(t))\otimes \Q$,  this group plays an important
role in topology in the work done by Casson and Gordon. 
\end{abstract}
\maketitle

 \section{$L$-groups and signatures} \label{applgroups}
Let $R$ be a ring with (possibly trivial) involution.
An $\eps$--hermitian ($\eps=\pm 1$)
form is a sesquilinear map $\t:V\times V\to R$ over a finitely generated free $R$--module
with the properties that $\t(ra,b)=\bar{r}\t(a,b), \t(a,rb)=\t(a,b)r$ and
$\t(a,b)=\eps \overline{\t(b,a)}$ for all $a,b\in V,r\in R$. It is called non--singular if the  map
$V\to
\hom(V,R), a\mapsto (b\mapsto \t(b,a))$ is an isomorphism. We denote by $L_0(R,\eps), \eps =\pm
1$, the Witt group of $\eps$--hermitian non--singular forms 
(cf. \cite{L93} and \cite{R98}). 
More precisely denote by $M$ the groupoid under direct sum  of $\eps$--hermitian non--singular forms. 
Let $\sim$ be the
equivalence relation generated by setting any form $(V,\t)$ to zero which has a submodule of half--rank
on which $\t$ vanishes. Then define $L_0(R,\eps):=M/\sim$, this is a group (cf. \cite{L93}) under the direct 
sum operation. We'll abbreviate
$L_0(R)$ for $L_0(R,+1)$.

%Any hermitian form $(V,\t)$ over $R$ can be represented,
%after choosing a basis for $V$,
%by a hermitian matrix $A$. This matrix is unique up to conjugation.
%Similarly an element in $L_0(R,\eps)$ can be represented by a matrix,
%but then the matrix is unique only up to conjugation and matrix
%cobordism.

Let  $F \subset \C$ be a subfield, closed under complex conjugation.
We will always equip the rings $F[t,t^{-1}], F(t)$
  with the involution given by the
complex involution on $F$ and $\bar{t}:=t^{-1}$.
For $\tau=[(V,\t)]\in L_0(F,\eps)$ we define
\[ \sign(\tau):=\dim(V^+)-\dim(V^-) \]
where $V^+$ (resp. $V^-$) denotes the maximal positive (resp. negative) subspace of $\sqrt{\eps}\t$.
 This number is independent of the choice
of representative
$(V,\t)$.

If $F \subset \C$ is a subfield such that all positive elements are squares, then
by Sylvester's theorem  
\[ \ba{rcl} \sign:{L}_0(F,\eps)&\to &\Z\\
(V,\t) &\mapsto&\sign(V,\t)\ea \]
is an isomorphism.
In particular  ${L}_0(\C,\pm 1)=L_0(\bar{\Q},\pm 1)\cong \Z$ via the signature map,
where we denote by $\bar{\Q}\subset \C$ the algebraic closure of $\Q$.
Since we are interested in studying to which degree signatures
determine forms we'll work in this paper with
$\ti{L}_0(R,\eps):=L_0(R,\eps)\otimes \Q$, i.e. we ignore the torsion part of $L_0(R,\eps)$.
Note that the above maps extend to an isomorphism $\sign:{L}_0(F,\eps)\otimes \Q \to \Q$.

Let $F$ be a Galois extension of $\Q$ with a (possibly trivial) involution.
Denote by $G(F)$ the
set of all $\Q$--linear embeddings
$F\to
\bar{\Q}$ preserving the involution.
For $\rho \in G(F)$ denote the induced maps
$\ti{L}_0(F,\eps)
\to \ti{L}_0(\bar{\Q},\eps), \ti{L}_0(F(t),\eps)\to \ti{L}_0(\bar{\Q}(t),\eps),\dots $
by $\rho$ as well.
Denote by  $G_0(F)\subset G(F)$ any subset such that for each
$\rho \in G(F)$ there exists precisely one
$\ti{\rho}\in G_0(F)$ with $\ti{\rho}=\rho$ or $\ti{\rho}=\bar{\rho}$.

Let $\tau=[(V,\t)]\otimes \frac{p}{q}\in \ti{L}_0(\bar{\Q}(t),\eps)$. If $z\in S^1$ is transcendental
then we consider $\C$ as a $\bar{\Q}(t)$ module via $f(t)\cdot w:=f(z)w$. Denote the
pairing 
\[ \ba{rcl} V\otimes_{\bar{\Q}(t)}\C\times V\otimes_{\bar{\Q}(t)}\C&\to&\C \\
(f_1(t)\otimes z_1,f_2(t)\otimes z_2)&\mapsto& \bar{z}_1 \t(f_1(t),f_2(t))|_{t=z}z_2    \ea \]
by $\t$ as well.
Then 
$\tau(z):=[(V\otimes_{\bar{\Q}(t)}\C,\t)]\otimes \frac{p}{q}\in \ti{L}_0(\C,\eps)$ is well--defined,
i.e. independent of the choice of $(V,\t)$.

The goal of this paper is to prove the following theorem.
\begin{theorem} \label{thmsignkt}
Let $F$ be a Galois extension of $\Q$, then for $\tau \in \ti{L}_0(F(t),\eps)$
we get
\[ \ba{rl} & \tau =0 \in \ti{L}_0(F(t),\eps)\\
\Leftrightarrow & \sign(\rho(\tau)(z))=0 \mbox{ for all }
 \rho \in G_0(F) \mbox{ and all transcendental } z\in S^1.\ea  \]
\end{theorem}

This result was stated by Litherland \cite[p. 358]{L84}, but there is no proof in the
literature.

{ \bf Acknowledgment.} I would like to thank Jerry Levine and Andrew Ranicki for helpful discussions.

%-------------------------
\section{Topological motivation}
The motivation for studying $\ti{L}_0(F(t))$ comes from knot theory. 
A knot $K\subset S^{n+2}$ is a smooth submanifold homeomorphic to $S^n$. A knot is called slice if it bounds
a smooth $(n+1)$--disk in $D^{n+3}$. For even dimensions Kervaire \cite{K65} showed that all knots are slice
and Levine \cite{L69} showed that in the case $n>1$ a knot is slice if and only if it is algebraically slice, i.e.
if its Seifert  form is metabolic. 

The case $n=1$ turned out to be much harder. A significant breakthrough was made by Casson and Gordon \cite{CG86} who found
examples of knots that are algebraically slice but not geometrically slice. They used a sliceness obstruction
which lies in $\ti{L}_0(F(t))$ where $F/\Q$ is a number field. We give a quick exposition of their results.

Let $K\subset S^3$ be an oriented knot. Denote 
by $M_K$ the result of zero framed surgery along $K$. 
There exists a canonical isomorphism $\eps:H_1(M_K)\to \Z$.
Let $k$ be   some prime power, denote by $M_k$ the $k$-fold cover corresponding to
$H_1(M_K)\to
\Z\to
\Z/k$. Denote by $TH_1(M_k)$ the $\Z$--torsion part of $H_1(M_k)$.

Let   $\chi:H_1(M_k) \to TH_1(M_k)  \to S^1$ be a character of order $m$.
Set
$F:=\Q(e^{2\pi i/m})$.
Since $\Omega_3(\Z \times \Z/m)=H_3(\Z \times \Z/m)$
is torsion (cf. \cite{CF64}) there exists $(V_k^4,\eps \times \chi)$ and some $r\in \N$ such
that
$\partial(V_k,\eps \times \chi)=r(M_k,\eps \times \chi)$. 
Let $\a :\Z \times S^1\to \mbox{Aut}(F(t))$ 
be the map given by $\a(n,z)(f(t)):=t^nzf(t)$.
Casson and Gordon showed that the $F(t)$--valued pairing on $H_2^{\a\circ (\eps \times \chi)}(V_k,F(t))$ is
non--singular and therefore defines an element $t(V_k)\in L_0(F(t))$. Denote by $t_0(V_k)$ the image of the
ordinary intersection pairing on $H_2(V_k)$ under the canonical map $L_0(\Z)\to L_0(F(t))$.

\begin{theorem}\cite{CG86}  \label{thmcgslice}
Let $K$ be a knot, $k,m$ prime powers and $\chi:H_1(L_k)\to S^1$ a character of order $m$.
\bn
\item 
 \[ \tau(K,\chi):=(t(V_k)-t_0(V_k))\otimes \frac{1}{r} \in \ti{L}_0(F_{\chi}(t))  \]
is a well-defined invariant of $(M_k,\eps \times \chi)$,
i.e. independent of the choice of $V_k$.
\item
  $H_1(M_k)=\Z \oplus TH_1(M_k)$.
\item
If furthermore $K$ is slice then there exists a subgroup  $Q \subset TH_1(M_k)$ with $|Q|^2=|TH_1(M_k)|$
  such that for any $\chi:TH_1(M_k)\to S^1$ of order $m$ vanishing on $Q$ we get $\tau(K,\chi)=0$.
\en
\end{theorem}

In   \cite{F03} the author shows that theorem \ref{thmcgslice}
can be completely reformulated in terms of eta invariants. This proof
relies in particular on the fact that elements in $\ti{L}_0(F(t))$ can be detected by signature invariants, i.e. 
the proof relies on theorem \ref{thmsignkt}.

\cmt{-----------------------------}
\section{Proof of the main theorem}
\subsection{The groups $\ti{L}_0(F), \ti{L}_0(F(t))$}

We quote a result from Ranicki \cite[p. 493]{R98}.
\begin{proposition} \label{prop91}
The following map is an isomorphism
\[ \ba{rcl} \ti{L}_0(F)
& \to & \bigoplus_{\rho \in G_0(F)} \ti{L}_0(\bar{\Q})  
  \\
\tau \otimes r & \mapsto & (\rho(\tau)\otimes r)_{\rho \in G_0(F)}  
\ea
\] 
%More precisely
%\[ L_0(F) \cong \Z^{|G_0(F)|} \oplus \mbox{8-torsion} \]
%In particular, $g=0 \in \ti{L}_0(F,\eps)$ if and only if
%$g=0 \in L_0(F,\eps)\otimes_{\Z} \Q$.
\end{proposition}

Consider the case where $F:=\Q[t]/q(t)$, $q(t)$ irreducible and  $q(t)=u q(t^{-1})$ for
some unit $u \in
\Q[t^{-1},t]$. Then there's an involution given by
$\bar{t}=t^{-1}$, which is non-trivial if $q(t) \ne t-1,t+1$. In this case the set
$G(F)$ corresponds canonically to the set of all roots of $q(t)$ lying in $S^1$ and $G_0(F)$
corresponds to all roots $z \in S^1$ of $q(t)$ with $\im(z) \geq 0$.

The goal of this section is to prove the following theorem.

\begin{theorem} \label{thmlfc}
Let $F$ be a Galois extension of $\Q$, then for $\tau \in \ti{L}_0(F(t),\eps)$
we get
\[  \tau =0 \in \ti{L}_0(F(t),\eps)
\Leftrightarrow \rho(\tau)=0 \in \ti{L}_0(\bar{\Q}(t),\eps) \mbox{ for all } \rho \in G_0(F)  \]
\end{theorem}

 To simplify the notation we'll only prove the case $\eps=1$.
We need some definitions and results from \cite[ch. 39C]{R98}.

\begin{definition}
Let $F$ be a field with a possibly trivial involution. Then define
$LAut_{fib}^0(F,\eps)$
to be the
Witt group of triples $(V,\theta,f)$ where $V$ is a
vector space over $F$, $\theta$ an $\eps$-hermitian form on $V$ and $f$ an isometry of $(V,\theta)$
 such that $(f-1)$ is an automorphism as well. Let 
$\ti{L}Aut_{fib}^0(F,\eps):=LAut_{fib}^0(F,\eps)\otimes \Q$.
\end{definition}

\begin{proposition} \label{prop93}  \cite[p. 533]{R98}
Let $F$ be a field with (possibly trivial) involution.
\bn
\item There exists a split exact sequence
\[
0 \to  \ti{L}_0(F[t,t^{-1}],\eps)  \to   \ti{L}_0(F(t),\eps) \to
\ti{L}Aut_{fib}^0(F,-\eps) \to  0\]
\item
Denote by $\overline{\mathfrak{M}}(F)$ the set of irreducible
monic polynomials $p(t)$ in $F[t]$ with
the added property that $\overline{p(t)}=up(t)$
for some unit $u \in F[t,t^{-1}]$ and
$\overline{\mathfrak{M}}^0(F):=\overline{\mathfrak{M}}(F)\setminus \{t-1\}$.
For $p(t)\in \overline{\mathfrak{M}}^0(F)$ define
\[ \ba{rrcl} r_{p(t)}:&\ti{L}Aut^0_{fib}(F,\eps)&\to&
\ti{L}_0(F[t,t^{-1}]/p(t),\eps)\\
 &(V,\t,f)\otimes r&\to& (\ker\{p(f):V \to V\},\ti{\t})\otimes r
\ea \]
where $\ti{\t}(a,b)=\sum_{i=0}^{deg(p)-1}\t(a,bt^i)t^{-i}$ and $t$ acts by $f$.
Then
\[\prod_{p(t) \in \overline{\mathfrak{M}}^0(F)}r_{p(t)} :  \ti{L}Aut_{fib}^0(F,\eps) \xrightarrow{\cong}
  \Oplus_{p(t)\in \overline{\mathfrak{M}}^0(F)} \ti{L}_0(F[t,t^{-1}]/p(t),\eps )\]
is an isomorphism and the inverse map is given by
\[\ba{rcl} \ti{L}_0(F[t,t^{-1}]/p(t),\eps )&\to&
 \ti{L}Aut_{fib}^0(F,\eps)\\
  (V,\theta)\otimes r&\mapsto&
 (V,\tr_{ (F[t,t^{-1}]/p(t))/F}\circ \theta,t)\otimes r \ea
 \]
\item The map \[ \ba{rcl} \ti{L}_0(F,\eps) &\to & \ti{L}_0(F[t,t^{-1}],\eps) \\
      (V,\theta)\otimes r &\mapsto &((V,\theta)\otimes_{F}F[t,t^{-1}])\otimes r \ea \]
is an isomorphism.
\en
\end{proposition}

There exists a commuting diagram of exact sequences ($G_0=G_0(F)$)
\[ \ba{rcccccccl}
0 &\to & \ti{L}_0(F[t,t^{-1}]) & \to &  \ti{L}_0(F(t)) &\to &
\ti{L}Aut_{fib}^0(F,-1) &\to & 0 \\[0.2cm]
\displaystyle && \downarrow \prod_{\rho \in G_0}\rho&& \downarrow \prod_{\rho \in
G_0}\rho&&\downarrow\prod_{\rho \in G_0}\rho \\[0.2cm]
\displaystyle 0 &\to & \bigoplus_{\rho \in G_0} \ti{L}_0(\bar{\Q}[t,t^{-1}]) &
\to &   \bigoplus_{\rho \in G_0} \ti{L}_0(\bar{\Q}(t)) &\to &\bigoplus_{\rho \in G_0}
 \ti{L}Aut_{fib}^0(\bar{\Q},-1) &\to & 0
\ea
 \]

From propositions \ref{prop91} and \ref{prop93}
it follows that the first vertical map is an injection.
Once we show that the last vertical map is an injection
as well it follows that the middle vertical map is an injection,
this will prove theorem \ref{thmlfc}.

For $p\in F[t,t^{-1}]$ irreducible we'll write $F_p:=F[t,t^{-1}]/p(t)$.
Note that there exists a canonical correspondence
\[ \ba{rcl}
\{ (\rho,z) | \rho \in G(F) \mbox{ and } z \in S^1 \mbox{ such that } \rho(p)(z)=0\}
&\leftrightarrow &G(F_p) \\
(\rho,z) &\mapsto &(\rho_{z}:\sum a_i t^i\to \rho(a_i)z^i) \ea \]
since $F/\Q$ is Galois.
Consider  
\[  \ba{rcl}
 \ti{L}Aut_{fib}^0(F,\eps) &\xrightarrow{\cong} &
 \displaystyle \Oplus_{p\in \overline{\mathfrak{M}}^0(F)}
\ti{L}_0(F_p,\eps )
 \hookrightarrow  \Oplus_{\rho \in G_0} \Oplus_{p\in
\overline{\mathfrak{M}}_0(F)}\Oplus_{\scriptsize{\ba{c} z \in S^1\sm \{1\} \\
\rho(p)(z)=0 \ea }}
 \ti{L}_0(\bar{\Q},\eps) \xrightarrow{\mu_{\rho,z}} \Q
\\
\downarrow\prod_{\rho \in
G_0}\rho \\
\displaystyle \bigoplus_{\rho \in G_0} \ti{L}Aut_{fib}^0(\bar{\Q},\eps) &\xrightarrow{\cong} &
\bigoplus_{\rho \in G_0} \Oplus_{z\in S^1\sm\{1\}}
\ti{L}_0(\bar{\Q},\eps) \xrightarrow{\s_{\rho,z}} \Q
\ea
 \]
where $\mu_{\rho,z}$ and $\sigma_{\rho,z}$ denotes the composition of
projection maps on the corresponding $\ti{L}_0(\bar{\Q},\eps)$ summand and taking
signatures. Note that $\mu_{\rho,z}$ is well-defined, since different $p(t)$'s have disjoint zero sets. Define
$\mu_{\rho,z}$ to be the zero map if $z$ is not a root for any $\rho(p(t))$.

\begin{proposition} \label{propsigmamuequal}
Let $p \in \overline{\mathfrak{M}}_0(F)$. For $(V,\theta) \in \ti{L}_0(F_p,\eps)$ we get
$\s_{\rho,z}(V,\theta)=\mu_{\rho,z}(V,\theta)$
for all $\rho
\in G_0(F)$ and $z \in S^1 \sm \{1\}$ such that $\rho(p)(z)=0$.
\end{proposition}

\begin{corollary}
The map
\[ \prod_{\rho \in G_0(F)}\rho: \ti{L}Aut_{fib}^0(F,\eps)
\to \bigoplus_{\rho \in G_0(F)} \ti{L}Aut_{fib}^0(\bar{\Q},\eps) \]
is an injection.
\end{corollary}

Note that this corollary concludes the proof of theorem \ref{thmlfc}.
We'll first prove the corollary.

\begin{proof}
The induced map
\[ \prod_{\rho,z\in S^1\sm \{1\}}\mu_{\rho,z}: \ti{L}_0(F_p,\eps) \to
\bigoplus_{\rho,z\in S^1 \sm \{1\}} \Q\]
is an injection. From the proposition it also follows that the induced map
\[ \prod_{\rho\in G_0(F),z\in S^1\sm \{1\}}\sigma_{\rho,z}: \ti{L}_0(F_p,\eps) \to \bigoplus_{\rho\in
G_0(F),z\in S^1}
\Q\] is an injection. Since different $p$'s have disjoint sets of zeros it follows that
\[ \prod_{\rho\in G_0(F),z\in S^1\sm \{1\}} \sigma_{\rho,z}: \Oplus_{p\in \overline{\mathfrak{M}}^0(F)}
 \ti{L}_0(F_p,\eps) \to \bigoplus_{\rho\in G_0(F),z\in
S^1} \Q\]
is an injection as well. But this implies that the intermediate map
\[ \prod_{\rho \in G_0(F)}\rho: \ti{L}Aut_{fib}^0(F,\eps)
\to \bigoplus_{\rho \in G_0(F)} \ti{L}Aut_{fib}^0(\bar{\Q},\eps) \]
is injective.
\end{proof}

Now we'll prove proposition  \ref{propsigmamuequal}.

\begin{proof}
Let $p\in \overline{\mathfrak{M}}_0(F)$. Denote the zeros of $\rho(p)(t)$ by $\a_1,\dots,\a_n$.
Pick a zero $\a$ of $\rho(p)$, we can assume that $\a=\a_1$.
Denote the induced embedding $F_p \to \bar{\Q}$ by $\rho_{\a}$.
Consider

\[ \ba{rccccccccl}
\ti{L}_0(\bar{\Q},\eps) &\hspace{-.15cm} \xleftarrow{r_{t-z}} \hspace{-.15cm}& \ti{L}Aut_{fib}^0(\bar{\Q},\eps)
& \hspace{-.15cm}\xleftarrow{\rho} \hspace{-.15cm}&   \ti{L}Aut_{fib}^0(F,\eps)&
\hspace{-.15cm}\leftarrow \hspace{-.15cm}& \ti{L}_0(F_p,\eps)
&\hspace{-.15cm} \xrightarrow{\rho_{\a}}\hspace{-.15cm} &\ti{L}_0(\bar{\Q},\eps)\\
\theta_l &\hspace{-.15cm}\mapstoleft\hspace{-.15cm} &\bar{\Q}\otimes_{F}(V,tr_{F_p/F}\circ
\theta,t)&\hspace{-.15cm} \mapstoleft \hspace{-.15cm}&(V,tr_{F_p/F}\circ
\theta,t)&\hspace{-.15cm}\mapstoleft\hspace{-.15cm}&(V,\theta)&\hspace{-.15cm}\mapsto\hspace{-.15cm} &
\theta_r
\ea  \]
We have to show that $\sign(\t_l)=\sign(\t_r)$.
Note that $\theta_r$ denotes the form
\[ \ba{rcl} \theta_r: V\otimes_{F_p}\bar{\Q} \times V\otimes_{F_p}\bar{\Q} & \to& \bar{\Q}\\
(v_1 \otimes_{F_p} z_1,v_2\otimes_{F_p} z_2) &\mapsto & \bar{z}_1\rho_{\a}(\theta(v_1,v_2))z_2 \ea
\]
where $F_p$ acts on $\bar{\Q}$ via $\rho_{\a}$.

Now we have to understand $\theta_l$.
In the following we'll view $\bar{\Q}$ as an $F$-module via $\rho$.
The form $\bar{\Q}\otimes_{F}(V,tr_{F_p/F}\circ \theta,t)$ is given by
\[ \ba{rcccccl}
V\otimes_{F}\bar{\Q} \times V\otimes_{F}\bar{\Q}   &\hspace{-.2cm}\xrightarrow{\tr_{F_p/F} \circ
\theta}\hspace{-.2cm}&
\bar{\Q}\otimes_{F}F\otimes_{F}\bar{\Q} &\hspace{-.2cm}\to \hspace{-.2cm}&
\bar{\Q}   \\
(v_1 \otimes_F z_1,v_2 \otimes_F z_2) &\hspace{-.2cm}\mapsto\hspace{-.2cm}&
\bar{z}_1\otimes_F \tr_{F_p/F}(\theta(v_1,v_2)) \otimes_F z_2  &\hspace{-.2cm}\mapsto\hspace{-.2cm}&
\bar{z}_1 \rho(\tr_{F_p/F}(\theta(v_1,v_2)))  z_2 \ea
\]

Denote by $\bar{\Q}_{p}$ the ring $F_p \otimes_F \bar{\Q} =\bar{\Q}[t,t^{-1}]/\rho(p(t))$.
It is easy to see that the map\[ \bar{\Q}_{p(t)} = F_p \otimes_F \bar{\Q} \xrightarrow{\tr_{F_p/F} \otimes_F \id} F
\otimes_F \bar{\Q} \to
\bar{\Q}\]
coincides with $\tr_{\bar{\Q}_{p}/\bar{\Q}}:\bar{\Q}_p \to \bar{\Q}$.
Therefore the form $\bar{\Q}\otimes_{F}(V,tr_{F_p/F}\circ \theta,t)$ is given by
\[ \ba{rcccccl}
V\otimes_{F}\bar{\Q} \times V\otimes_{F}\bar{\Q} &\to & \bar{\Q}_{p} &\xrightarrow{\tr_{\bar{\Q}_p/\bar{\Q}}}
&\bar{\Q} \\
(v_1 \otimes_F z_1,v_2 \otimes_F z_2)&\mapsto&  \theta(v_1,v_2) \otimes_F \bar{z}_1z_2 &\mapsto&
\tr_{\bar{\Q}_p/\bar{\Q}}(\theta(v_1,v_2) \otimes_F \bar{z}_1z_2) \ea
\]
We can write $V \otimes_{F} \bar{\Q}=V_1 \oplus \dots \oplus V_n$ where
$V_i:=\ker \{(t-\a_i):V\otimes_{F}\bar{\Q} \to V\otimes_{F} \bar{\Q} \}$ since the minimal polynomial of $t$ is
$p(t)=\prod_{i=1}^n (t-\a_i)$. Then $\t_l$ is given by restricting
the above form to $V_1$.

We can decompose the $\bar{\Q}[t]$-module $\bar{\Q}_{p}=\bar{\Q}[t]/\rho(p(t))$ as follows
\[   \bar{\Q}_p = \bigoplus_{i=1}^n \ker\{(t-\a_i):\bar{\Q}_{p}\to \bar{\Q}_{p}\}
    = \Oplus_{i=1}^n \bar{\Q}_i
        \]
where $\bar{\Q}_i:=\ker\{(t-\a_i):\bar{\Q}_p\to \bar{\Q}_p\}$.
Note that $\dim_{\bar{\Q}}(\bar{\Q}_i)=1$. Consider the  map $\mu_{\a_i}: \bar{\Q}_i \to \bar{\Q}$
given by $p(t) \mapsto p(\alpha_i)$. These maps define isomorphisms of $\bar{\Q}$-algebras.
Then the trace function is given by
\[ \ba{rrcl} \tr_{\bar{\Q}_{p}/\bar{\Q}}:& \oplus_{i=1}^n \bar{\Q}_i=\bar{\Q}_{p} &\to &\bar{\Q} \\
      &  (z_1,\dots,z_n) &\mapsto & \sum_{i=1}^n \mu_{i}(z_i) \ea \]
since $\tr_{\bar{\Q}_{p}/\bar{\Q}}=\tr_{(\Oplus_{i=1}^n \bar{\Q}_i)/\bar{\Q}}=\sum_{i=1}^n
\tr_{\bar{\Q}_i/\bar{\Q}}$.
The form  $\theta \otimes_F \bar{\Q}:V \otimes_F \bar{\Q} \times V \otimes_F \bar{\Q} \to \bar{\Q}_p$ restricts to
a form $V_1 \times V_1 \to \bar{\Q}_1$ and $\t_l$ is given by
$V_1 \times V_1 \to \bar{\Q}_1 \xrightarrow{\tr}\bar{\Q}$.

We can now compute $\theta_l$. Let
$\sum_{j=1}^{s_1} v_{1j} \otimes_F z_{1j}, \sum_{l=1}^{s_2} v_{1l}
\otimes_F z_{1l} \in V_1$, then
\[ \ba{rcl}
\t_l(\sum_{j=1}^{s_1} v_{1j} \otimes_F z_{1j}, \sum_{l=1}^{s_2} v_{1l}
\otimes_F z_{1l})
&=& \tr_{\bar{\Q}_p/\bar{\Q}}(\sum_{j=1}^{s_1}\sum_{l=1}^{s_2} \theta(v_{1j},v_{2l})
\otimes_F \bar{z}_{1j}z_{2l})=
\\ &=&
\mu_{\a}(\sum_{j=1}^{s_1}\sum_{l=1}^{s_2} \rho(\theta(v_{1j},v_{2l}))  \bar{z}_{1j}z_{2l})
=\\
&=&\sum_{j=1}^{s_1}\sum_{l=1}^{s_2} \rho_{\a}(\theta(v_{1j},v_{2l}))  \bar{z}_{1j}z_{2l}
\ea \]

Consider the following sequence of canonical isomorphisms:
\[ \ba{rccccccl}
V_1=\ker\{(t-\a):V\otimes_F \bar{\Q}\hspace{-.1cm}\to\hspace{-.1cm} V\otimes_{F}\bar{\Q}\}
&\hspace{-.2cm}\cong \hspace{-.2cm}&(V\otimes_{F}\bar{\Q})\otimes_{\bar{\Q}[t]}\bar{\Q}
&\hspace{-.2cm}\cong\hspace{-.2cm} &V\otimes_{F[t]}\bar{\Q} &\hspace{-.2cm}\cong\hspace{-.2cm} & V\otimes_{F_p}
\bar{\Q}  \ea \] The resulting isomorphism  is  given by
\[ \ba{rccccccl}
V_1=\ker\{(t-\a):V\otimes_F \bar{\Q}\to V\otimes_{F}\bar{\Q}\}  &\cong &
V\otimes_{F_p} \bar{\Q}
 \\
\sum_{j=1}^s v_j \otimes_{F} z_j &\mapsto & \sum_{j=1}^s v_j\otimes_{F_p} z_j
\ea \]
It now follows immediately that the forms $\theta_l,\theta_r$ are isomorphic.
\end{proof}

\cmt{-----------------------------}
\subsection{The group $\ti{L}_0(\bar{\Q}(t))$}
We need to quote one more fact.

\begin{proposition} \label{prop94}  \cite[p. 533]{R98}
Let $F$ be a field. 
The splitting $ \ti{L}Aut_{fib}^0(F,-\eps) \to \ti{L}_0(F(t),\eps)$ in the exact sequence
\[ 
0 \to  \ti{L}_0(F[t,t^{-1}],\eps)  \to   \ti{L}_0(F(t),\eps) \to
\ti{L}Aut_{fib}^0(F,-\eps) \to  0\]
 is given by
\[
(V,\theta,f) \mapsto  (V\otimes_{F}F(t),
(v,w)\to (1-t^{-1})\t((1-f)^{-1}v,w)+ 
\eps (1-t)\overline{\t((1-f)^{-1}w,v)})
 \]
\end{proposition}

\begin{theorem} \label{thmlct}
Let $\tau \in \ti{L}_0(\bar{\Q}(t))$, then
\[   \tau =0 \in \ti{L}_0(\bar{\Q}(t))
\Leftrightarrow \tau(z)=0 \in \ti{L}_0(\bar{\Q}) \mbox{ for all transcendental } z.
 \]
\end{theorem}

 Combining theorems \ref{thmlfc} and \ref{thmlct} we now get a proof for theorem \ref{thmsignkt}.

\begin{proof}
Proposition  \ref{prop93}, part (1) and (3), shows that there exists an isomorphism
\[ \ti{L}_0(\bar{\Q}) \oplus \ti{L}Aut_{fib}^0(\bar{\Q},-1) \to \ti{L}_0(\bar{\Q}(t)) \]
Let $Z:=S^1\sm \{1\} \cap \bar{\Q}$, then $\overline{\mathfrak{M}}_0(\bar{\Q})=\{t-z | z \in Z\}$.
Using that  $\ti{L}_0(\bar{\Q})\cong \Q$ via the signature and 
using part (2) of proposition
 \ref{prop93} we get 
isomorphisms
\[ \ba{rcccl} \bigoplus_{z \in Z} \Q & \xrightarrow{\cong} &
\Oplus_{z \in Z} \ti{L}_0(\bar{\Q}[t,t^{-1}]/(t-z),-1) &\xrightarrow{\cong}& \ti{L}Aut_{fib}^0(\bar{\Q},-1)\\
     (r_{z})_{z \in Z}
&\mapsto &  \bigoplus (\bar{\Q}[t,t^{-1}]/(t-z),i)\otimes r_z
 &\mapsto & \bigoplus (\bar{\Q},i,z)\otimes r_z \ea \]

The isomorphisms above and proposition \ref{prop94}  show that in $\ti{L}_0(\bar{\Q}(t))$
the form $\tau$ is equivalent to
\[  (\bar{\Q}(t),1)\otimes r_0 \oplus \bigoplus_{j=1}^s
(\bar{\Q}(t),i(1-t^{-1})(1-\bar{z}_j)^{-1}+i(1-t)(1-{z_j})^{-1})\otimes r_j
\]
where $z_j \in S^1\sm \{1\},j=1,\dots,s$ are distinct and $r_0,\dots,r_s\in \Q$.
Note that $\tau=0$ if and only if
$r_0=r_1=\dots=r_s=0$.

We can assume that $r_i\in \N$ for all $i$, and hence restrict ourselves to forms
in $L_0(\bar{\Q}(t))$. Then the matrix
\[  A(t):= (1)\otimes r_0 \oplus \bigoplus_{j=1}^s
(i(1-t^{-1})(1-\bar{z}_j)^{-1}+i(1-t)(1-{z_j})^{-1})\otimes r_j
\]
represents $\tau$.
The signature function $z \mapsto \sign(A(z))$ is locally constant,
its only jumps are when $\det(A(t))=0$,
i.e. when
\[ (1-t^{-1})(1-\bar{z}_j)^{-1}+(1-t)(1-{z_j})^{-1}=0 \mbox{ for some } j \]
i.e. when $t=\frac{1-\bar{z}_j}{1-{z_j}} \in S^1$. It is clear that
$\sign(A(1))=r_0$ and that the jump of the signature function at $=\frac{1-\bar{z}_j}{1-{z_j}} \in
S^1$ is $2r_j$. The theorem follows now easily since the transcendental numbers in $S^1$ are dense.
\end{proof}

\end{document}